\documentclass[12pt]{article}
\usepackage[usenames]{color}
\usepackage{graphicx, subfigure}
\usepackage{amsmath, amsthm, amssymb}
\usepackage{amsfonts}
\usepackage{fullpage}
\usepackage{ifthen}
\usepackage{url}
\usepackage[sort&compress]{natbib}
\usepackage{multirow}
\usepackage{bm}
\usepackage{tikz-qtree}

\usepackage[linesnumbered,algoruled,boxed,lined,commentsnumbered]{algorithm2e}

\makeatletter
\def\@biblabel#1{}
\makeatother

\theoremstyle{plain}
\theoremstyle{definition}
\theoremstyle{remark}

\title{Why must every data scientist be a Platonist}
\author{Collin Joseph Van Cuyk$^*$\qquad Liang Hong\footnote{Department of Mathematical Sciences, The University of Texas at Dallas, 800 West Campbell Road, Richardson, TX 75080, USA.}}

\usepackage{tikz}
\usetikzlibrary{automata,arrows,positioning,calc}
\usepackage{listings}
\usepackage[labelformat=empty]{caption}
\usepackage{subfig}
\usepackage{caption}
\begin{document}

\maketitle

\begin{abstract}
Data scientists are not mathematicians, but they make heavy use of mathematics in their daily work.  While mathematicians can study a mathematical object which is inaccessible to our five senses, data scientists  must deal with real-world data which are observable to us.  This fine line suggests that a data scientist's philosophical position on mathematics might have a nontrivial impact on their work.  By examining how different philosophical views of mathematics affect the interpretation of the basic model assumption in data science, we arrive at the conclusion that a data scientist, who uses modern probabilistic and statistical tools, must be a Platonist.

\smallskip

\end{abstract}

\section{Introduction}

Mathematics is an indispensable tool for every data scientist.  While most data scientists barely examine their own philosophical positions on mathematics,  various philosophical outlooks challenge the basic assumptions made by data scientists.  In the current literature, there are three dominating schools of thought concerning mathematical philosophy:
(i) Platonism, (ii) Finitism, and (iii) Formalism.  The names of these three schools vary slightly in the literature; see, for example,  Benacerraf and Putnam (1983), George and Vellenman (2002), Henle (1992),  Kleene (1950),  and Kunen (1980, 2011, 2012). In this article, we will follow the terminologies of Chapter~III of Kunen (2012).  Specifically,  we will consider a Platonist to be the one who considers all of mathematics literally real.  Apparently,  there is a good portion of modern mathematics which is not observable in our real world, such as many well-founded sets.  A Platonist believes this part of mathematics exists in a world that is inaccessible to our five senses.  Next, let a finitist be one who works with only those mathematical entities which are finite.   Finitists reject the notion of a definite infinity, although they will allow for the concept of a potential infinite,  but clarify that this process will never actually be completable. In short, a finitist believes only in finitary mathematics. Finally,  we take a formalist be one who accepts as real only those things which a finitist accepts, but recognizes the intellectual intrigue of infinite mathematics. That is, they only recognize the things accepted within finitism as meaningful, but they consider any infinite mathematics as an intellectual game.  While none of these three philosophical positions is considered to be the ``correct'' position,  we will show that each has a nontrivial impact on the basic model assuption in data science.  In particular, we will draw the conclusion that an honest data scientist must be a Platonist.  Throughout, by a data scientist, we mean a person whose profession is to analyze data using usual tools in modern probability and statistics. 

\section{The basic model assumption in data science}

The most fundamental assumption of data science is that there is an underlying unobservable probability law that governs all observable data.  Depending on the problem at hand, the data scientist chooses his assumption about this law.  Examples include but are not limited to a density model, a linear model,  or a time series model.  To illustrate our point, it suffices to consider the parametric density model.  Most data science textbooks describe such a model as follows.  Assume that data $X_1, \ldots, X_n$ are independently and identically distributed (iid) according to a true (but unobservable) probability density function $f^\star$.  Since $f^\star$ is unobservable,  a data scientist posits a parametric model 
\[\mathcal{M}=\{f_{\theta}\mid \theta\in \Theta \},\]
where $\Theta$, called the parameter space, is a subset of a finite-dimensional Euclidean space.  For example, if we consider a gamma model with shape parameter $\alpha>0$ and rate parameter $\lambda>0$ where both $\alpha$ and $\lambda$ are unknown,  then
\[ \mathcal{M}=\{f_{(\alpha, \lambda)}\mid (\alpha, \lambda)\in \mathbb{R}_+\times \mathbb{R}_+\}, \]
where $\mathbb{R}_+=(0, \infty)$ and 
\[f_{(\alpha, \lambda)}(x)=\frac{\lambda^\alpha}{\Gamma(\alpha)}x^{\alpha-1}  e^{-\lambda x}, x>0,\]
and $\Gamma$ is the gamma function.

If a data scientist is a Platonist, then this gamma model $\mathcal{M}$ makes perfect sense,  every ingredient used in the model description is well-defined, and $\mathcal{M}$  is essentially an uncountable family of gamma distributions.  That is, a Platonist's interpretation of this model is the same as what is in a data science textbook (e.g.  Lehmann 1983 and Schervish 1995).

However,  to a finitist, the above model description cannot be interpreted as the textbook intends.  What's even worse,  this model makes no sense to a finitist. First, since a finitist rejects any definite infinity,  the parameter space of the gamma model, i.e. $\mathbb{R}_+\times \mathbb{R}_+$,  cannot be an uncountable set.  A finitist can only allow the parameter space $\Theta$ to be a finite set.  This evidently differs from what is described in a data science textbook; see, for example,  Lehmann (1983) and Schervish (1995).  A moment of reflection reveals that there is a more serious issue a finitist will take exception to.  In modern probability theory,  a probability density function $f$  is an almost surely nonnegative Borel-meausrable function such that $\int_{\infty}^{-\infty}f(t)dt=1$, and 
\[F(x)=\int_{-\infty}^x f(t)dt,\]
where $F$ is the corresponding cumulative distribution function.  Since a definite integral is the limit of a Riemann sum which is not an acceptable concept  to a finitist, the notion of a probability density function does not make sense to a finitist.  Therefore, the description of a gamma model using gamma densities make no sense to a finitist, either.  Similarly,  a finitist will reject the nice relationship between of a probability density function and its corresponding distribution function $F'(x)=f(x)$. Therefore, a data scientist cannot be a finitist. 

How about a formalist? Can a data scientist be a formalst? A formalist is a ``diplomatic'' person; he takes a position to avoid being confrontational with either a Platonist or  a finitist.  Though his position is beyond reproach in mathematics,  the story is different when it comes to data science.  A formalist only believes finitary mathematics to be real. Therefore, the above gamma model makes no sense to him either.  Can he simply consider this model as an intellectual model (which might be meaningless in the real-world) as he always does for the infinitary mathematics? The answer is negative.  Since a data scientist  deals with real-world data which are observable to us,  his job is to propose a ``best'' model to describe these real-world data. Therefore, an honest data scientist cannot posit a model he does not believe exists.

\section{What about Bayesian analysis?}

In the previous section, we evidently took the frequentist approach. However, if a data scientist options for the Bayesian approach,  he will face the same issue.  To see that,  we recall that the parametric Bayesian model treats data $X_1, \ldots, X_n$ to be iid conditioning on the parameter $\theta$ and puts a prior distribution $\Pi_0$ on $\theta$. That is,  a parametric Bayesian model is specified hierarchically as 
\[\theta\sim \Pi_0 \text{ and $X_1, \ldots, X_n\mid \theta \overset{iid}{\sim} f_{\theta}$},\]
where $f_{\theta}$ is the probability density function on the $X$-space indexed by $\theta\in \Theta$.  Then the prior distribution $\Pi_0$ is updated to the posterior distribution $\Pi_n$ via the Bayes formula:
\[\Pi_n (A)=\frac{\int_A L(\theta) d\Pi_0(\theta)}{\int_{\Theta} L(\theta)d\Pi_0(\theta)},\quad A\subset \Theta,\]
where $L(\theta)=\prod_{i=1}^n f(X_i)$ is the likelihood function and $A$ is a $\Pi_0$-measurable set.  For a detailed treatment of Bayesian analysis, we refer to Chapter 1 of Schervish (1995) and Chapter 4 of Berger (2006).  Because integration and probability density functions (hence the likelihood function) are both involved in a Bayesian model, the same argument we used in the previous section shows that a data scientist can be  neither a finitist nor a formalist.  Thus a data scientist must be a Platonist. 

\section{Concluding remarks}
If a data scientist hopes to be consistent with their own philosophy, there is little they, as a finitist, could accomplish. Most if not all practical applications of statistical inference fail to function within absolutely defined domains or purely discrete, purely finite models. It is not our purpose to suggest which philosophy has more merit, however for all practical or applied purposes the data scientist cannot restrict themselves so significantly as required by finitism. What of the formalist? Believing to be real only mathematics accepted by the finitist, the formalist will struggle with the same contradictions and problems between their work and their philosophy. The formalist could not ethically present findings which they considered to be only “a game”. We are forced, then, to conclude that the modern data scientist must hold the philosophy of Platonism if they wish to complete their work satisfied that it is both real and meaningful.

%\section*{Acknowledgments}

%The authors thank the Editor, Associate Editor, and the two anonymous referees for their comments and suggestions, which led to %significant improvements in this article.

\section*{References}

\begin{description}

\item{} Berger, J.O.~(2006).  \emph{Statistical Decision Theory and Bayesian Analysis}. Springer: New York.

\item{} Benacerraf,  P.~and Putnam,  H.~(1983). \emph{Philosophy of mathematics,} Second Edition. Cambridge University Press: Cambridge, UK. 

%\item{} Enderton, H.B.~(1977). \emph{Elements of Set Theory}. Academic Press: San Diego, CA.

\item{} George, A.~and Velleman, D. J.~(2002).  \emph{Philosophies of Mathematics.}. Blackwell Publishers Inc. 

\item{} Henle, J.M.~(1991). The happy formalist. \emph{Mathematical Intelligencer}~13(1), 12--18.

%\item{} Jech, T.~(2003). \emph{Set Theory}.  Springer: New York.

\item{} Kleene, S.C. ~(1950).  \emph{Introduction to MetaMathematics}. North-Holland: Amsterdam. 

\item{} Kunen, K.~(1980).  \emph{Set Theory: An Introduction to Independence Proofs}. Elsevier: Amsterdam

\item{} Kunen, K.~(2011). \emph{Set Theory}. College Publications: London, UK.

\item{} Kunen, K.~(2012). \emph{The Foundations of Mathematics}, Revised Edition. College Publications: London, UK.

\item{} Lehmann, E.L.~(1983). \emph{Theory of Point Estimation}. Wiley: New York.

%\item{} Smullyan, R.M.~and Fitting, M.~(2010).  \emph{Set Theory and the Continuum Problem}. Dover Publications, Inc: New York.

%\item{} Russell, B.~(1919). \emph{Introduction to Mathematical Philosophy},  George Allen $\&$ Unwin.

%\item{} Shoenfield, J.R. (1967). \emph{Mathematical Logic}. CRC Press: Boca Raton, London, New York.

\item{} Schervish, M.J.~(1995). \emph{Theory of Statistics}. Springer: New York.

\end{description}

\end{document}